\documentclass[12pt,a4paper]{article}
\usepackage[width=18cm,height=23cm]{geometry}
\usepackage{amsmath,amssymb,amsfonts,amsthm,enumerate}
\providecommand{\keywords}[1]{\textbf{\textit{Keywords: }}#1}

\usepackage{graphicx}
\usepackage{tikz}
\usetikzlibrary{shapes,backgrounds}
\usepackage{tikz-cd}
\usetikzlibrary{%
  matrix,%
  calc,%
  arrows%
}

\usepackage[all]{xy}
\usepackage[utf8]{inputenc}
\usepackage[T1]{fontenc}
\usepackage{mathrsfs,amsbsy,bm,mathtools}
\usepackage{indent first}

\newcommand{\cd}{\mathcal{D}}

\newcommand{\bc}{\mathbb{C}}

\newcommand{\bz}{\mathbb{Z}}

\newtheorem{Thm}{Theorem}[section]%
\newtheorem{Prop}[Thm]{Proposition}%
\theoremstyle{definition}
\newtheorem{Def}[Thm]{Definition}%
\newtheorem{Rmk}[Thm]{Remark}%
\newtheorem{Quest}[Thm]{Question}%

\title{\bf Affine Homogeneous Surfaces with Hessian rank 2
\\
and Algebras of Differential Invariants}
\author{Zhangchi Chen, Jo\"el Merker}
\date{\today}

\begin{document}
\maketitle

\abstract{Consider a graphed holomorphic surface $u=F(x,y)$ in $\bc^3_{x,y,u}$ under the action of the affine transformation group $A(3)$. In 1999, Eastwood and Ezhov obtained a list of homogeneous models by determining possible tangential vector fields. Inspired by Olver's recurrence formulas, we study the algebra of $A(3)$ differential invariants of surfaces. We obtain necessary conditions for homogeneity of algebraic nature. Solving these conditions, we organise homogeneous models in inequivalent branches.}


\section{Introduction}
In continuation with~{\cite{Olver-2007, Chen-Merker-2020,
Arnaldsson-Valiquette-2020}}, we study the algebras of
differential invariants of surfaces in $3$-dimensional space, under
the affine $A(3)$\,\,---\,\,or special affine
$SA(3)$\,\,---\,\,transformation group. To fix the context, we work
over $\mathbb{C}$, and we study surfaces $S^2 \subset \mathbb{C}^3$ under $A(3)$ whose
Hessian has (maximal) rank $2$.  We investigate ramifications of the
Lie-Tresse theorem, following the theory of
Fels-Olver~{\cite{Fels-Olver-1998, Fels-Olver-1999}}.

Mainly, we explore what the powerful {\sl recurrence relations}
provide, firstly to determine the structures of the concerned 
algebras of differential invariants, and
secondly to determine 
all $A(3)$-homogeneous nondegenerate surfaces $S^2 
\subset \mathbb{C}^3$. The homogeneous classification
was terminated in~{\cite{Abdalla-Dillen-Vrancken-1997,
Doubrov-Komrakov-Rabinovich-1996,
Eastwood-Ezhov-1999}}, without regard to algebras
of differential invariants.

Interestingly, we really need to know the {\em explicit expressions}
of certain key (relative or absolute) differential invariants which
create {\sl bifurcation branches}.

Thus, consider a holomorphic local surface $S^2$ in
$\mathbb{C}^3 \ni (x,y,u)$ graphed as
\begin{equation}
\label{u-F-x-y-initial}
u
=
F(x,y)
=
\sum\limits_{j+k\geqslant 0}\,
F_{j,k}\,
\frac{x^j}{j!}\,\frac{y^k}{k!}.
\end{equation}
Using translations of $A(3)$, we may assume $F(0,0) = 0$, so $j + k
\geqslant 1$.  The goal is to normalize the power series coefficients
$F_{j,k}$ using the $GL(3)$ action.

The Hessian determinant $\big\vert \begin{smallmatrix} F_{xx} & F_{xy}
\\ F_{yx} & F_{yy}
\end{smallmatrix} \big\vert$ 
is a $GL(3)$-relative invariant, and we assume it is nowhere
vanishing.  After elementary $GL(3)$ transformations done in
Section~{\ref{sect-normalization}}, we can {\em pre}normalize $u = F$ to
\[
u
=
x\,y
+
G_{3,0}\,
\frac{x^3}{6}
+
G_{0,3}\,
\frac{y^3}{6}
+
\sum\limits_{j+k\geqslant 4}\,
G_{j,k}\,
\frac{x^j}{j!}\,\frac{y^k}{k!},
\]
where all the $G_{j,k} = 
G_{j,k} 
(F_{{\scriptscriptstyle{\bullet}}, {\scriptscriptstyle{\bullet}}})$
express in terms of the $F_{l,m}$ with $l + m \leqslant j+k$.
On a computer, we store these (long) expressions.

The stabilizer of such a {\em pre}normal form consists of
bi-dilations $(x, y, u) \longmapsto \big( \mu x, \lambda y, 
\mu \lambda u \big)$, with $\lambda, \mu \in \mathbb{C}^\ast$,
and of the swap $x \longleftrightarrow y$. 
Consequently, $G_{3,0}$ and $G_{0,3}$, and even all the
higher order $G_{j,k}$, are relative 
invariants\footnote{\,We would like to mention that the product
$G_{3,0}\, G_{0,3}$ is a nonzero multiple of the so-called
{\sl Pick invariant}, whose identical vanishing
characterizes {\em ruled} surfaces. The numerator of the product
$G_{3,0}\, G_{0,3}$ is shown in the beginning of Section~{\ref{sect-normalization}}. The numerator of $G_{3,0}$ is shown in \eqref{explicit-30}.}.
 
Admitting Lie's principle that any (relative) invariant can be assumed
either $\equiv 0$ or $\neq 0$ after restriction to some open subset,
$G_{3,0}$ and $G_{0,3}$ create $3$ branches, up to $x
\longleftrightarrow y$.

\smallskip\noindent${\bf B_1}$:\,
$G_{3,0} \neq 0$ and $G_{0,3} \neq 0$.

\smallskip\noindent{$\bf B_2$}:\,
$G_{3,0} \neq 0$ and $G_{0,3} \equiv 0$.

\smallskip\noindent{$\bf B_3$}:\,
$G_{3,0} \equiv 0$ and $G_{0,3} \equiv 0$.

\medskip

Abbreviating `${\sf root}$' to denote the Hessian rank $2$ assumption
$F_{1,1}^2 - F_{2,0} F_{0,2} \neq 0$, here is the complete branching
diagram to which the next five statements will refer.
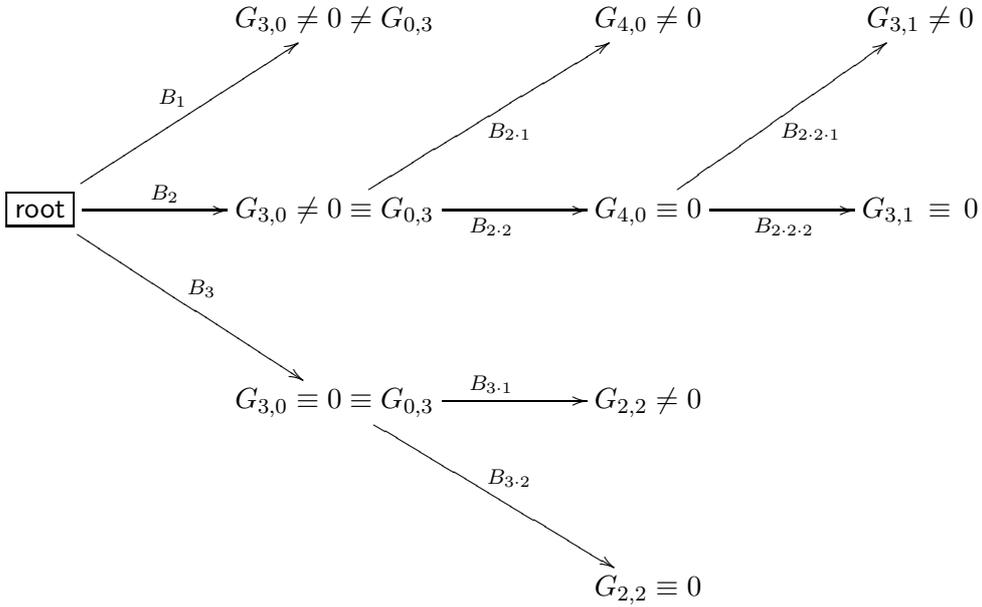
\[
\xymatrix{
&&
G_{3,0}\neq0\neq G_{0,3}
&&
G_{4,0}\neq0
&&
G_{3,1}\neq0
&&
\\
&&
&&
&&
&&
\\
\ar[uurr]^{B_1}
\ar[rr]^{B_2\,\,\,\,\,\,\,\,\,\,\,\,}
\boxed{\sf root}
\ar[ddrr]^{B_3}
&&
\ar[uurr]_{B_{2\cdot 1}}
\ar[rr]_{B_{2\cdot 2}}
G_{3,0}\neq0\equiv G_{0,3}
&&
\ar[uurr]_{B_{2\cdot2\cdot 1}}
\ar[rr]_{B_{2\cdot2\cdot2}}
G_{4,0}\equiv 0
&&
G_{3,1}\,\equiv\,0
&&
\\
&&
&&
&&
&&
\\
&&
G_{3,0}\equiv 0\equiv G_{0,3}
\ar[rr]^{B_{3\cdot 1}}
\ar[rrdd]^{B_{3\cdot 2}}
&&
G_{2,2}\neq 0
&&
&&
\\
&&
&&
\\
&&
&&
G_{2,2}\equiv 0
}
\]
In this tree, any two surfaces landing in one of the six different
terminal branches are always $A(3)$-{\em in}equivalent.

\begin{Thm}
In the first branch $B_1$ where $G_{3,0}
(F_{{\scriptscriptstyle{\bullet}}, {\scriptscriptstyle{\bullet}}})
\neq 0 \neq G_{0,3}
(F_{{\scriptscriptstyle{\bullet}}, {\scriptscriptstyle{\bullet}}})$,
the following hold.

\smallskip\noindent{\bf (1)}\,
The graphed equation
normalizes as
\[
u
=
x\,y
+
\frac{x^3}{6}
+
\frac{y^3}{6}
+
\sum_{j+k\geqslant 4}\,
I_{j,k}
(F_{{\scriptscriptstyle{\bullet}},{\scriptscriptstyle{\bullet}}})\,
\frac{x^j}{j!}\,
\frac{y^k}{k!},
\]
where all $I_{j,k}$ are differential invariants,
up to the swap $x \longleftrightarrow y$ and a discrete group
\[
\mathcal{G}_0:=\big\{x'=\omega^j\,x, y'=\omega^{-j}\,y, u'=u~|~j=0,1,2\big\}
\]
where $\omega:=e^{2\pi i/3}$, a cube root of unity.

\smallskip\noindent{\bf (2)}\,
The algebra of differential invariants is generated by
$I_{4,0}$, $I_{3,1}$, $I_{1,3}$, $I_{0,4}$ and all their invariant derivatives
$\mathcal{D}_1^{\alpha_1} \mathcal{D}_2^{\alpha_2} 
({}_{{}^{{}^{\scriptscriptstyle{\bullet\!}}}})$,
with $\alpha_1, \alpha_2 \in \mathbb{N}$. In particular, $I_{2,2}$ can be solved
\[
I_{2, 2} = \tfrac{8}{9}\,I_{4, 0}\,I_{0, 4}-\tfrac{1}{9}\,I_{1, 3}\,I_{3, 1}+\tfrac{2}{9}\,I_{4, 0}\,I_{3, 1}-\tfrac{1}{36}\,\cd_xI_{3, 1}+\tfrac{1}{36}\,\cd_yI_{4, 0}.
\]

\smallskip\noindent{\bf (3)}\,
The moduli space of all possible homogeneous models 
is described, in the space of coefficients
$\big( I_{4,0}, I_{3,1}, I_{2,2}, I_{1,3}, I_{0,4} \big)
\in \mathbb{C}^5$, by the complex algebraic variety of dimension $2$ defined
by
\[
\aligned
(E1) \ \ 0&=8\,I_{0, 4}\,I_{4, 0}-I_{1, 3}\,I_{3, 1}+2\,I_{3, 1}\,I_{4, 0}-9\,I_{2, 2},\\
(E2) \ \ 0&= 2\,I_{0, 4}\,I_{1, 3}+8\,I_{0, 4}\,I_{4, 0}-I_{1, 3}\,I_{3, 1}-9\,I_{2, 2},\\
(E3) \ \ 0&= 4\,I_{0, 4}\,I_{3, 1}-I_{1, 3}\,I_{2, 2}-4\,I_{2, 2}\,I_{4, 0}+2\,I_{3, 1}^2+9\,I_{1, 3}+18\,I_{4, 0},\\
(E4) \ \ 0&= 4\,I_{0, 4}\,I_{2, 2}-2\,I_{1, 3}^2-4\,I_{1, 3}\,I_{4, 0}+I_{2, 2}\,I_{3, 1}-18\,I_{0, 4}-9\,I_{3, 1}.
\endaligned
\]
\end{Thm}

Precisely, there is a one-to-one correspondence between
$A(3)$-equivalence {\em classes} of homogeneous surfaces $S^2 \subset
\mathbb{C}^3$ in branch $B_1$ and points $\big( I_{4,0}, I_{3,1}, I_{2,2},
I_{1,3}, I_{0,4} \big) \in \mathbb{C}^5$ satisfying $(E1)$, $(E2)$, $(E3)$,
$(E4)$, modulo the swap and $\mathcal{G}_0$.  In Section~{\ref{sect-nonvanishing-Pick}}, we resolve these equations and reobtain,
without overlap, models $N1$, $N2$, $N3$, $N4$
of~{\cite{Eastwood-Ezhov-1999}}.

It is elementary to verify that any affine vector field
which is tangent to the surface is a linear combination of the two
independent ones
\[
\aligned
e_1
&
:=
\big(
1
-
\tfrac{1}{2}\,I_{2,2}\,u
+
\tfrac{1}{4}\,u
-
\tfrac{1}{3}\,I_{1,3}\,x
-
\tfrac{2}{3}\,I_{4,0}\,x
\big)\,
\partial_x
+
\big(
-\tfrac{1}{2}\,I_{3,1}\,u
-
\tfrac{2}{3}\,I_{1,3}\,y
-
\tfrac{1}{3}\,I_{4,0}\,y
-
\tfrac{1}{2}\,x
\big)\,\partial_y
\\
&
\ \ \ \ \ \ \ \ \ \ \ \ \ \ \ \ \ \ \ \ \ \ \ \ \ \ \ \ \ \ \ \ \ \ 
\ \ \ \ \ \ \ \ \ \ \ \ \ \ \ \ \ \ \ \ \ \ \ \ \ \ \ \ \ \ \ \ \ \ 
\ \ \ \ \ \ \ \
+
\big(
-I_{1,3}\,u
-
I_{4,0}\,u
+
y
\big)\,\partial_u,
\\
e_2
&
:=
\big(
-\tfrac{1}{2}\,I_{1,3}\,u
-
\tfrac{1}{2}\,y
-
\tfrac{2}{3}\,I_{3,1}\,x
-
\tfrac{1}{3}\,I_{0,4}\,x
\big)\,\partial_x
+
\big(
1
-
\tfrac{1}{2}\,I_{2,2}\,u
+
\tfrac{1}{4}\,u
-
\tfrac{1}{3}\,I_{3,1}\,y
-
\tfrac{2}{3}\,I_{0,4}\,y
\big)\,\partial_y
\\
&
\ \ \ \ \ \ \ \ \ \ \ \ \ \ \ \ \ \ \ \ \ \ \ \ \ \ \ \ \ \ \ \ \ \ 
\ \ \ \ \ \ \ \ \ \ \ \ \ \ \ \ \ \ \ \ \ \ \ \ \ \ \ \ \ \ \ \ \ \ 
\ \ \ \ \ \ \ \ \ 
+
\big(
-I_{0,4}\,u-I_{3,1}\,u+x
\big)\,\partial_u.
\endaligned
\]
Moreover, computing the Lie bracket $[e_1, e_2]$ and
subtracting appropriate linear combinations
of $e_1$ and $e_2$ to get a vector field vanishing
at the origin, this pair
of vector fields {\em does constitute} a $2D$ Lie algebra
with the uniquely defined Lie bracket:
\[
[e_1,e_2]
=
\big(
-\tfrac{2}{3}\,I_{3,1}
-
\tfrac{1}{3}\,I_{0,4}
\big)\,e_1
+
\big(
\tfrac{1}{3}\,I_{4,0}
+
\tfrac{2}{3}\,I_{1,3}
\big)\,
e_2,
\]
{\em if and only if} equations $(E1)$, $(E2)$, $(E3)$,
$(E4)$ hold.

\smallskip

All the other branches 
$B_{2\cdot 1}$, $B_{2\cdot 2\cdot 1}$,
$B_{2\cdot 2\cdot 2}$, $B_{3\cdot 1}$, $B_{3\cdot 2}$
have $0 \equiv G_{0,3} 
(F_{{\scriptscriptstyle{\bullet}}, {\scriptscriptstyle{\bullet}}})$.
Similarly as for the study of {\em parabolic} surfaces
(constant Hessian rank $1$) achieved in~{\cite{Chen-Merker-2020}},
it is necessary to {\em insert} this differential relation
and its consequences into the 
power series normalizations and into
all recurrence relations as well. 
In Section~{\ref{sect-v03-nv40}}, we introduce the 
relevant notion of {\sl subjets}.
As a matter of fact, $G_{0,3} 
(F_{{\scriptscriptstyle{\bullet}}, {\scriptscriptstyle{\bullet}}})
\equiv 0$ can be solved as $F_{y^3} = R_{0,3}
({}_{{}^{{}^{\scriptscriptstyle{\bullet\!}}}})$,
with some complicated remainder, whence
all derivatives $F_{x^j y^k}$ with $k \geqslant 3$ are
{\sl dependent}. 

\begin{Thm}
In the second
branch $B_{2\cdot 1}$ where $G_{3,0} \neq 0 \equiv G_{0,3}$ 
and $G_{4,0} \neq 0$, 
the following holds.

\smallskip\noindent{\bf (1)}\,
The graphed equation normalizes as
\[
u
=
x\,y
+
\frac{x^3}{6}
+
\frac{y^3}{6}
+
\frac{x^{24}}{24}
+
I_{3,1}\frac{x^3y}{6}
+
I_{2,2}\frac{x^2y^2}{4}
+
I_{1,3}\frac{xy^3}{6}
+
I_{0,4}\frac{y^4}{24}
+
\sum_{j+k\geqslant 5}\,
I_{j,k}\,
(F_{{\scriptscriptstyle{\bullet}},{\scriptscriptstyle{\bullet}}})\,
\frac{x^j}{j!}\,\frac{y^k}{k!},
\]
where all $I_{j,k}$ are differential invariants.

\smallskip\noindent{\bf (2)}\,
The algebra of differential invariants is generated by $I_{3,1}$, $I_{2,2}$, $I_{5,0}$ and all their invariant derivatives
$\mathcal{D}_1^{\alpha_1} \mathcal{D}_2^{\alpha_2}
({}_{{}^{{}^{\scriptscriptstyle{\bullet\!}}}})$,
with $\alpha_1, \alpha_2 \in \mathbb{N}$. In particular, $I_{4,1}$ can be solved
\[
I_{4, 1} = -8\,I_{3, 1}^2+2\,I_{5, 0}\,I_{3, 1}+\cd_x\,I_{3, 1}+\tfrac{7}{2}\,I_{2, 2}-2\,I_{3, 1}.
\]

\smallskip\noindent{\bf (3)}\,
The moduli space of all possible homogeneous models
is exactly described, in the space
$\mathbb{C}^3 \ni \big( I_{3,1}, I_{2,2}, I_{5,0}, I_{4,1}, I_{3,2} \big)$,
by the complex-algebraic variety
of dimension $1$ defined by the $4 + 3$ equations
\[
\aligned
(E41) \ \ 0& =I_{4,1}+8I_{3,1}^2-\tfrac{7}{2}I_{2,2}+2I_{3,1}-2I_{5,0}I_{3,1},
\\
(E42) \ \ 0&=4I_{3,1}I_{2,2}+2I_{3,1}^2-2I_{4,1}I_{3,1}+I_{3,2},
\\
(E43) \ \ 0&=12I_{3,1}I_{2,2}-3I_{5,0}I_{2,2}+4I_{2,2}+I_{3,2},
\\
(E44) \ \ 0&=6I_{2,2}^2+4I_{3,1}I_{2,2}-3I_{4,1}I_{2,2},
\endaligned
\]
\[
\aligned
(F51) \ \ 0&=24\,I_{3, 1}^2\,I_{5, 0}-2\,I_{5, 0}^2\,I_{3, 1}-\tfrac{15}{2}\,I_{2, 2}\,I_{5, 0}+7\,I_{5, 0}\,I_{3, 1}
\\
& \ \ \ +\tfrac{21}{2}\,I_{2, 2}-64\,I_{3, 1}^3+36\,I_{2, 2}\,I_{3, 1}-40\,I_{3, 1}^2-6\,I_{3, 1},
\\
(F52) \ \ 0&=30\,I_{2, 2}\,I_{3, 1}+72\,I_{2, 2}\,I_{3, 1}^2-18\,I_{2, 2}\,I_{5, 0}\,I_{3, 1}-\tfrac{63}{4}\,I_{2, 2}^2
\\
& \ \ \ +56\,I_{3, 1}^3-14\,I_{3, 1}^2\,I_{5, 0}+12\,I_{3, 1}^2+64\,I_{3, 1}^4-32\,I_{3, 1}^3\,I_{5, 0}+4\,I_{5, 0}^2\,I_{3, 1}^2,
\\
(F53) \ \ 0&=-I_{3, 1}\,(-16\,I_{3, 1}^2+4\,I_{3, 1}\,I_{5, 0}+3\,I_{2, 2}-6\,I_{3, 1})\,(-32\,I_{3, 1}^2+8\,I_{3, 1}\,I_{5, 0}+6\,I_{2, 2}-13\,I_{3, 1}).
\endaligned
\]
\end{Thm}

The Lie symmetry algebra is always of dimension $2$,
generated by
\[
\aligned
e_1
&
:=
\big[
1
+
(1-I_{5,0}+4I_{3,1})\,x
-
\tfrac{1}{2}I_{2,2}u
\big]\,
\partial_x
+
\big[
-\tfrac{1}{2}x
+
(3-2I_{5,0}+8I_{3,1})\,y
-
\tfrac{1}{2}I_{3,1}u
\big]\,
\partial_y
\\
&
\ \ \ \ \ \ \ \ \ \ \ \ \ \ \ \ \ \ \ \ \ \ \ \ \ \
\ \ \ \ \ \ \ \ \ \ \ \ \ \ \ \ \ \ \ \ \ \ \ \ \ \
\ \ \ \ \ \ \ \ \ \ \ \ \ \ \ \ \ \ \ \ \ \ \ \ \ \
+
\big[
y
+
(4-3I_{5,0}+12I_{3,1})\,u
\big]\,
\partial_u,
\\
e_2
&
:=
\big[
(I_{3,1}-I_{4,1}+2I_{2,2})\,x
\big]\,
\partial_x
+
\big[
1
+
(3I_{3,1}+4I_{2,2}-2I_{4,1})\,y
-
\tfrac{1}{2}I_{2,2}u
\big]\,
\partial_y
\\
&
\ \ \ \ \ \ \ \ \ \ \ \ \ \ \ \ \ \ \ \ \ \ \ \ \ \
\ \ \ \ \ \ \ \ \ \ \ \ \ \ \ \ \ \ \ \ \ \ \ \ \ \
\ \ \ \ \ \ \ \ \ \ \ \ \ \ \ \ \ \ \ \ \ \ \ \ \ \
+
\big[
x
+
(4I_{3,1}+6I_{2,2}-3I_{4,1})\,u
\big]\,
\partial_u,
\endaligned
\]
having Lie bracket
\[
[e_1,e_2]
=
\big(
I_{3,1}-I_{4,1}+2I_{2,2}
\big)\,e_1
+
\big(
-3+2I_{5,0}-8I_{3,1}
\big)\,
e_2,
\]
if and only if the above $7$ equations are satisfied.
We solve these equations and recover 
models $N5$ and $N6$ of~{\cite{Eastwood-Ezhov-1999}}.

\smallskip

In the next branch $B_{2 \cdot 2}$,
two differential relations exist, $G_{0,3}
(F_{{\scriptscriptstyle{\bullet}}, {\scriptscriptstyle{\bullet}}})
\equiv 0$
and   $G_{4,0}
(F_{{\scriptscriptstyle{\bullet}}, {\scriptscriptstyle{\bullet}}})
\equiv 0$. These two PDEs can be solved as
$F_{y^3} = R_{0,3} 
({}_{{}^{{}^{\scriptscriptstyle{\bullet\!}}}})$
and 
$F_{x^4} = R_{4,0} 
({}_{{}^{{}^{\scriptscriptstyle{\bullet\!}}}})$, with
complicated but explicit right-hand sides.
The power series reads:
\[
u
=
x\,y
+
\frac{x^3}{6}
+
G_{3,1}\,\frac{x^3y}{6}
+
G_{2,2}\,\frac{x^2y^2}{4}
+
G_{3,2}\,\frac{x^3y^2}{12}
+
{\sf dependent}\,{\sf remainder}.
\]
Only $3$ independent coefficients remain, and by analyzing the three
compatibility conditions
\[
\aligned
D_x^4 (R_{0,3}) &= D_y^3 (R_{4,0}),\\
D_x^5 (R_{0,3}) &= D_x\,D_y^3 (R_{4,0}),\\
D_x^6 (R_{0,3}) &= D_x^2\,D_y^3 (R_{4,0}),\\
\endaligned
\]
where $D_x$ and $D_y$ are the two total differentiation operators, we find
$G_{2,2} = 0$ and $G_{3,2} = 0$.  It remains only $G_{3,1}$.  But
there is still one degree of freedom $(x,y,u) \longmapsto (\mu x,
\mu^2 y, \mu^3 u)$ with $\mu \in \mathbb{C}^\ast$.  Then $G_{3,1}$ is
a relative differential invariant, which causes the two branches
$B_{2\cdot 2 \cdot 1}$ and $B_{2 \cdot 2 \cdot 2}$.

\begin{Prop}
In the third branch $B_{2 \cdot 2 \cdot 1}$ where
$G_{3,0} \neq 0 \equiv G_{0,3} \equiv G_{4,0}$ and
$G_{3,1} \neq 0$, the normal form
\[
u
=
x\,y
+
\frac{x^3}{6}
+
\frac{x^3y}{6}
+
\frac{x^5}{30}
+
\sum_{j+k\geqslant 6}\,
I_{j,k}
(F_{{\scriptscriptstyle{\bullet}},{\scriptscriptstyle{\bullet}}})\,
\frac{x^j}{j!}\,
\frac{y^k}{k!},
\]
has all its coefficients uniquely determined and
is automatically homogeneous, with $2$-dimensional 
affine Lie symmetry algebra generated by
\[
e_1
\,:=\,
\partial_x
+
\big(-\tfrac{1}{2}u-\tfrac{1}{2}x)\,\partial_y
+
y\,\partial_u,
\ \ \ \ \ \ \ \ \ \ \ \ \ \ \ \ \ \ \ \
e_2
:=
(1+y)\,\partial_y
+
(x+u)\,\partial_u,
\]
having Lie bracket $[e_1, e_2] = 0$. A closed form is
\[
u
\,=\,
(1+y)\,\sqrt{2}\,
\tan\,
\big(
\tfrac{x}{\sqrt{2}}
\big)
-
x.
\]
\end{Prop}

In particular, there is no way of getting any `algebra'
of differential invariants, because
all $I_{j,k}$ are constant!
This is model $N7$ of~{\cite{Eastwood-Ezhov-1999}},
while model $N8$ is recovered by

\begin{Prop}
In the fourth branch $B_{2\cdot 2\cdot 2}$ where $G_{3,0} \neq 0
\equiv G_{0,3} \equiv G_{4,0} \equiv G_{3,1}$, the unique
normal form is Cayley's cubic
\[
u
\,=\,
x\,y
+
\frac{x^3}{6},
\]
with $3D$ affine symmetries
\[
e_1
\,:=\,
\partial_x-\tfrac{1}{2}\,x\,\partial_y+y\,\partial_u,
\ \ \ \ \ \ \ \ \ \ \ \ \ \ \ \ \ \ \ \
e_2
\,:=\,
\partial_y+x\,\partial_u,
\ \ \ \ \ \ \ \ \ \ \ \ \ \ \ \ \ \ \ \
e_3
\,:=\,
x\,\partial_x
+
2y\,\partial_y
+
3u\,\partial_u,
\]
having solvable Lie structure
$[e_1, e_3] = e_1$, $[e_2, e_3] = 2 e_2$.
\end{Prop}

Next, in the branch $B_3$,
when $G_{3,0} 
(F_{{\scriptscriptstyle{\bullet}},{\scriptscriptstyle{\bullet}}})
\equiv 0 \equiv G_{0,3} 
(F_{{\scriptscriptstyle{\bullet}},{\scriptscriptstyle{\bullet}}})$,
two PDEs of the form $F_{x^3} = R_{3,0}
({}_{{}^{{}^{\scriptscriptstyle{\bullet\!}}}})$ 
and  $F_{y^3} = R_{0,3}
({}_{{}^{{}^{\scriptscriptstyle{\bullet\!}}}})$ 
are satisfied
by $F(x,y)$, hence all $F_{x^jy^k}$ with either
$j \geqslant 3$ or $k \geqslant 3$ are dependent. 
Only $G_{2,2}$ remains unnormalized, and it is a relative
invariant, since 
the remaining freedom is
$(x,y,u) \longmapsto (\mu x, \lambda y, \mu \lambda u)$,
with $\mu, \lambda \in \mathbb{C}^\ast$.
Then $G_{2,2}$ creates the last two branches
$B_{3\cdot 1}$ and $B_{3 \cdot 2}$,
which we gather in a single statement.

\begin{Prop}
{\bf (a)}\,
In the fifth branch $B_{3 \cdot 1}$ where
$G_{3,0} \equiv G_{0,3} \equiv 0 \neq G_{2,2}$,  
there exists a single surface
\[
u
\,=\,
x\,y
+
\frac{x^2y^2}{4}
+
\sum_{j+k\geqslant 5}\,
I_{j,k}\,
\frac{x^j}{j!}\,
\frac{y^k}{k!},
\]
with uniquely determined $I_{j,k} \in \mathbb{C}$, which is
automatically homogeneous, having $3D$ symmetries
\[
e_1
\,:=\,
-\,x\,\partial_x+y\,\partial_y,
\ \ \ \ \ \ \ \ \ \ \ \ \ \ \ \ \ \ \ \
e_2
\,:=\,
(u-2)\,\partial_x-2y\,\partial_u,
\ \ \ \ \ \ \ \ \ \ \ \ \ \ \ \ \ \ \ \
e_3
\,:=\,
(u-2)\,\partial_y-2x\,\partial_u,
\]
with structure $\cong \mathfrak{sl}(2,\mathbb{C})$
given by $[e_1,e_2] = e_2$, 
$[e_1,e_3] = -\,e_3$, $[e_2,e_3] = -\,2\,e_1$.
A closed form is
\[
u
\,=\,
2
-
2\,\sqrt{1-xy}.
\]

\smallskip\noindent{\bf (b)}\,
In the sixth, last branch $B_{3 \cdot 2}$ where
$G_{3,0} \equiv G_{0,3} \equiv G_{2,2} \equiv 0$, 
the normal form is the basic quadric
\[
u
\,=\,
x\,y,
\]
having $4D$ symmetries
\[
e_1
\,:=\,
-\,x\,\partial_x+y\,\partial_y,
\ \ \ \ \ \ \ \ \ \ \ \ \ \ \ \ \ \ \ \
e_2
\,:=\,
x\,\partial_x+u\,\partial_u,
\ \ \ \ \ \ \ \ \ \ \ \ \ \ \ \ \ \ \ \
e_3
\,:=\,
\partial_x+y\,\partial_u,
\ \ \ \ \ \ \ \ \ \ \ \ \ \ \ \ \ \ \ \
e_4
\,:=\,
\partial_y+x\,\partial_u,
\]
with solvable Lie structure
$[e_1,e_3] = e_3$, 
$[e_1,e_4] = 
-\,e_4$,
$[e_2,e_3] = -\,e_3$.
\end{Prop}

In February 2020, at IHES, we came to the following

\begin{Quest}
{\sl 
Can one reconstitute the branching tree
of differential invariants in 
Cartan's classification~{\cite{Cartan-1932}} of homogeneous
real hypersurfaces $M^3 \subset \mathbb{C}^2$?}
\end{Quest}

Beyond, here is a more demanding question, which
requires in principle to go up to order $6$.

\begin{Quest}
{\sl Can one determine the branching tree
of differential invariants related to 
the classification(s) of
multiply transitive
real hypersurfaces $M^5 \subset \mathbb{C}^3$ 
due to Loboda~{\cite{Loboda-2020}}, 
and to Doubrov-Medvedev-The~{\cite{Doubrov-Medvedev-The-2020}}?}
\end{Quest}

\smallskip\noindent{\bf (1)}\,
Explain bifurcations caused by some specific (relative) 
differential invariants.

\smallskip\noindent{\bf (2)}\,
Set up appropriate recurrence relations taking account of
ambient subjets, and determine minimal generators.

\smallskip\noindent{\bf (3)}\,
Determine moduli spaces of homogeneous models,
solve equations, and classify sharply. 

\medskip\noindent{\bf Acknowledgments.}
The realization of this research work in differential invariants has
received generous financial support from the scientific grant 2018/29/B/ST1/02583 originating from the Polish National Science Center (NCN).

In the context of Cartan's method of equivalence, we learned from
Pawe{\l} Nurowski the naturality of branching trees of (relative)
differential invariants while ``{\sl hunting for}'' homogeneous models.

In March 2019, we benefited from Boris Doubrov's
visit in Orsay University. 
The second author also acknowledges 
enlightening zoom exchanges
with Dennis The and Boris Doubrov. 

\section{Normalization, relative invariants and branchings}\label{sect-normalization}
There is an order 3 relative invariant called {\sl Pick invariant}. If it is non-zero we can normalize it to 1. If it is zero by homogeneity we assume it is constant 0. The numerator is
\[
\aligned
P:=& \  \ \ 6 F_{yyy} F_{xy} F_{xyy} F_{xx}^2 - 9 F_{yy} F_{xyy}^2 F_{xx}^2 - F_{yyy}^2 F_{xx}^3 - 
 12 F_{yyy} F_{xy}^2 F_{xx} F_{xxy} + 18 F_{yy} F_{xy} F_{xyy} F_{xx} F_{xxy}\\
& + 
 6 F_{yy} F_{yyy} F_{xx}^2 F_{xxy} - 9 F_{yy}^2 F_{xx} F_{xxy}^2 + 8 F_{yyy} F_{xy}^3 F_{xxx} - 
 12 F_{yy} F_{xy}^2 F_{xyy} F_{xxx} - 6 F_{yy} F_{yyy} F_{xy} F_{xx} F_{xxx}\\
& + 6 F_{yy}^2 F_{xyy} F_{xx} F_{xxx} + 
 6 F_{yy}^2 F_{xy} F_{xxy} F_{xxx} - F_{yy}^3 F_{xxx}^2
\endaligned
\]

We will see Pick when we normalize the Taylor coefficients.

\subsection{First loop} Start from a general Taylor expansion of a graphed surface in $\bc^3$
\[
u=F(x,y)=\sum\limits_{j+k\geqslant 0}\frac{F_{j,k}}{j!k!}(x-x_0)^j(y-y_0)^k.
\]
After certain elementary $A(3)$ action we may assume
\begin{equation}\label{u-O-2}
u=O(2)=\frac{F_{2,0}}{2}x^2+F_{1,1}xy+\frac{F_{0,2}}{2}y^2+\sum\limits_{j+k\geqslant 3}\frac{F_{j,k}}{j!k!}x^jy^k.
\end{equation}

First, we look at a surface close to the form $u=x^2-y^2+O(3)$, i.e. $(F_{2,0},F_{1,1},F_{0,2})$ in a neighborhood of $(2,0,-2)$. We want to find an affine transformation
\[
x'=ax+by, \ \ y'=cx+dy, \ \ u'=u
\]
close to the identity, sending $u=O(2)$ as in~\eqref{u-O-2} to $u'=x'^2-y'^2+O(3)$. We may proceed as follows
\[
\aligned
u&=
\frac{F_{2,0}}{2}\Big(x+\frac{F_{1,1}}{F_{2,0}}y\Big)^2+\frac{F_{0,2}F_{2,0}-F_{1,1}^2}{2F_{2,0}}y^2\\
&=\Big(\underbrace{\sqrt{\frac{F_{2,0}}{2}}x+\frac{F_{1,1}}{\sqrt{2F_{2,0}}}y}_{x'}\Big)^2-\Big(\underbrace{\sqrt{\frac{F_{1,1}^2-F_{0,2}F_{2,0}}{2F_{2,0}}}y}_{y'}\Big)^2.
\endaligned
\]
The transformation
\[
x'=\sqrt{\frac{F_{2,0}}{2}}x+\frac{F_{1,1}}{\sqrt{2F_{2,0}}}y, \ \ y'=\sqrt{\frac{F_{1,1}^2-F_{0,2}F_{2,0}}{2F_{2,0}}}y, \ \ u'=u
\]
is well defined for $(F_{2,0},F_{1,1},F_{0,2})$ in a neighborhood of $(2,0,-2)$, and tends to the identity when $(F_{2,0},F_{1,1},F_{0,2})$ tends to $(2,0,-2)$.

Next, we look at a surface close to the form $u=xy+O(3)$, i.e. $(F_{2,0},F_{1,1},F_{0,2})$ in a neighborhood of $(0,1,0)$. Again we want to find an affine transformation close to the identity, sending $u=O(2)$ to $u'=x'y'+O(3)$. This can be done after the change of coordinates $(x,y)=(s+t,s-t)$ and $(x',y')=(s'+t',s'-t')$. Our original surface becomes
\[
\aligned
u&=\sum\limits_{j+k\geqslant 2}\frac{F_{j,k}}{j!k!}(s+t)^j(s-t)^k\\
&=\frac{F_{2,0}+2F_{1,1}+F_{0,2}}{2} s^2+ (F_{2,0}-F_{0,2}) st +\frac{F_{2,0}-2F_{1,1}+F_{0,2}}{2}t^2+O(3)
\endaligned
\]
while our target surface is $u'=s'^2-t'^2+O(3)$. Here the new coefficients $(F_{2,0}+2F_{1,1}+F_{0,2},F_{2,0}-F_{0,2},F_{2,0}-2F_{1,1}+F_{0,2})$ are in a neighborhood of $(2,0,-2)$. We plug them in the transformation above. We conclude that the transformation
\[
s'=\sqrt{\frac{F_{2,0}+2F_{1,1}+F_{0,2}}{2}} s+\frac{F_{2,0}-F_{0,2}}{\sqrt{2}\sqrt{F_{2,0}+2F_{1,1}+F_{0,2}}} t, \ \ t'= \sqrt{\frac{2 (F_{1,1}^2 - F_{2,0}F_{0,2})}{F_{2,0} + 2 F_{1,1} + F_{0,2}}} t,  \ \ u'=u
\]
does the job. Hence in coordinates $(x,y)$ and $(x',y')$, the transformation is
\[
\aligned
x'&=\frac{F_{2,0}+F_{1,1}+\sqrt{F_{1,1}^2-F_{2,0}F_{0,2}}}{\sqrt{2}\sqrt{F_{2,0}+2F_{1,1}+F_{0,2}}} x +\frac{F_{1,1}+F_{0,2}-\sqrt{F_{1,1}^2-F_{2,0}F_{0,2}}}{\sqrt{2}\sqrt{F_{2,0}+2F_{1,1}+F_{0,2}}} y,\\
y'&=\frac{F_{2,0}+F_{1,1}-\sqrt{F_{1,1}^2-F_{2,0}F_{0,2}}}{\sqrt{2}\sqrt{F_{2,0}+2F_{1,1}+F_{0,2}}} x+\frac{F_{1,1}+F_{0,2}+\sqrt{F_{1,1}^2-F_{2,0}F_{0,2}}}{\sqrt{2}\sqrt{F_{2,0}+2F_{1,1}+F_{0,2}}} y, \ \ u'=u
\endaligned
\]
normalizes the surface to a new graph
\[
u'=G(x',y')=x'y'+\sum\limits_{j+k\geqslant 3}\frac{G_{j,k}}{j!k!}x'^jy'^k,
\]
where we can solve $G_{j,k}$ in terms of $F_{j,k}$. In particular
\[
\aligned
G_{3,0}=\frac{P_{36}}{2\sqrt{2} (F_{2,0}+ 2 F_{1,1} + F_{0,2})^{\frac{3}{2}} (F_{1,1}^2 - F_{2,0} F_{0,2})^{\frac{3}{2}}},
\endaligned
\]
where $P_{36}$ is a polynomial in $\bz[F_{j,k},\sqrt{F_{1,1}^2-F_{2,0}F_{0,1}}]$ with 36 monomials.

\subsection{Second loop}
Now we look at a general surface of the form
\[
u=F(x,y)=xy+O(3)=xy+\sum\limits_{j+k\geqslant 3}\frac{F_{j,k}}{j!k!}x^jy^k.
\]
The stabilizer, unfortunately, is not connected. To stabilise the form $u=xy+O(3)$, one uses either
\[
\mathcal{G}_1:=\{x'=\mu x+ku, y'=\lambda y+lu, u'=\mu \lambda u~|~\mu,\lambda\in\bc^*\},
\]
a dimension 4 subgroup, or
\[
\mathcal{G}_2:=\{x'=\mu y+ku, y'=\lambda x+lu, u'=\mu \lambda u~|~\mu,\lambda\in\bc^*\},
\]
a dimension 4 coset of $\mathcal{G}_1$ by switching $x$ and $y$.

We first study the effect of $\mathcal{G}_1$ acting on the third order Taylor coefficients. After a transformation
\[
x'=\mu x+ku, \ \ y'=\lambda y+lu, \ \ u'=\mu\lambda u,
\]
we obtain a new graph
\[
u'=G(x',y')=x'y'+\sum\limits_{j+k\geqslant 3}\frac{G_{j,k}}{j!k!}x'^jy'^k.
\]
Now the fundamental equation
\[
G\big(\mu x+kF(x,y),\lambda y+lF(x,y)\big)=\mu\lambda F(x,y)
\]
for $(x,y)$ in a neighborhood of $(0,0)$, lead us to solve
\begin{equation}\label{G-30-21-12-03}
\aligned
G_{3,0}&= \frac{\lambda F_{3,0}}{\mu^2},\\
G_{2,1}&=\frac{-2 l + F_{2,1} \lambda}{\lambda \mu},\\
G_{1,2}&=\frac{-2 k + F_{1,2} \mu}{\lambda \mu},\\
G_{0,3}&=\frac{ \mu F_{0,3}}{\lambda^2}.
\endaligned
\end{equation}
We see that $G_{3,0}$ and $G_{0,3}$ are relative invariants under $\mathcal{G}_1$, while $G_{2,1}$ and $G_{1,2}$ can be normalized to 0 by a unique choice of $l$ and $k$ (depends on $\lambda$ and $\mu$).

If we switch $x$ and $y$, we switch $G_{3,0}$ and $G_{0,3}$, $G_{2,1}$ and $G_{1,2}$. Consequently, $G_{3,0}G_{0,3}$ is a $\mathcal{G}_1$ and $\mathcal{G}_2$ relative invariant, hence an $A(3)$ relative invariant. In terms of the original Taylor coefficients $u=\sum\limits_{j+k\geqslant 2}\frac{F_{j,k}}{j!k!}x^jy^k$, we rediscover the numerator of the Pick invariant shown in the beginning of this section
\[
G_{3,0}G_{0,3}=\frac{-P}{8(F_{1,1}^2-F_{2,0}F_{0,2})^3},
\]

\begin{Rmk} The numerators of $G_{3,0}$ and of $G_{0,3}$ are polynomials having 36 monomials and are non-factorizable. But their product factorizes by $P$ (13 monomials) and $(F_{2,0}+2F_{1,1}+F_{0,2})^3$. The reason is that the polynomial ring $\bc[F_{j,k},\sqrt{F_{2,0}+2F_{1,1}+F_{0,2}},\sqrt{F_{1,1}^2-F_{2,0}F_{0,2}}]$ is not a UFD. The author thanks Z. Jiang (U. Michigan) for giving the following example:

Let $X:=F_{2,0}$, $Y:=F_{1,1}$, $S:=\sqrt{F_{2,0}+2F_{1,1}+F_{0,2}}$ and $T:=\sqrt{F_{1,1}^2-F_{2,0}F_{0,2}}$. Then $S^2+X T^2-(X+Y)^2=0$. Thus $(S+X+Y)(S-X-Y)=-XT^2$. That ring is not a UFD.
\end{Rmk}

\begin{Rmk} Indeed the polynomial ring $\bc[F_{j,k},\sqrt{F_{1,1}^2-F_{2,0}F_{0,2}}]$ is not a UFD.

Let $X:=\sqrt{F_{1,1}^2-F_{2,0}F_{0,2}}$. Then $F_{1,1}^2-X^2=F_{2,0}F_{0,2}$. Thus $(F_{1,1}-X)(F_{1,1}+X)=F_{2,0}F_{0,2}$. That ring is not a UFD.
\end{Rmk}

As a consequence of the remark
Hence we get 3 branches.
\begin{itemize}
\item {$\bf B_1$:} both $G_{3,0}$ and $G_{0,3}$ are non-zero. This is Pick non-vanishing branch.
\item {$\bf B_2$:} only $G_{3,0}$ is non-zero while $G_{0,3}\equiv0$, so Pick vanishes.
\item {$\bf B_3$:} both $G_{3,0}$ and $G_{0,3}$ are identically zero, so Pick vanishes as well.
\end{itemize}

In $B_1$, we normalize $G_{3,0}$ and $G_{0,3}$ to $1$ by a unique choice of $\mu$ and $\lambda$. We conclude that every non-degenerate surface with non-vanishing Pick is $A(3)$-equivalent to a graph
\[
u=xy+\frac{x^3}{6}+\frac{y^3}{6}+O(4).
\]
The stabilizer is a discrete group
\[
\mathcal{G}_0:=\big\{x'=\omega^j\,x, y'=\omega^{-j}\,y, u'=u~|~j=0,1,2\big\}
\]
where $\omega:=e^{2\pi i/3}$, a cube root of unity. We will study this branch further in the next section.

In $B_2$, we normalize $G_{3,0}$ to $1$ by a unique choice of $\lambda$ depending on $\mu$. Thus the group element
\[
l=\frac{F_{2,1}}{2 F_{3,0}},\ \  k=\frac{F_{1,2}}{2}, \ \  \lambda=\frac{1}{F_{3,0}}, \ \ \mu=1,
\]
sends $u=xy+\sum\limits_{j+k\geqslant 3}\frac{F_{j,k}}{j!k!}x^jy^k$ to $u'=x'y'+\frac{x'^3}{6}+O(4)$. Again we can solve all $G_{j,k}$ in terms of $F_{j,k}$. For $j+k=4$
\[
\aligned
G_{0,3}&=F_{0,3}F_{3,0}^2 \ \ \text{(in Branches 2 and 3, both sides are 0)},\\
G_{4,0}&=-\frac{2 F_{2,1} F_{3,0} - F_{4,0}}{F_{3,0}},\\
G_{3,1}&=
\frac{-3 F_{2,1}^2 - 4 F_{1,2} F_{3,0} + 2 F_{3,1}}{2},\\
G_{2,2}&=-3 F_{1,2} F_{2,1} F_{3,0} + 
  F_{2,2} F_{3,0},\\
G_{1,3}&= 
\frac{-3 F_{1,2}^2 F_{3,0}^2 + 2 F_{1,3} F_{3,0}^2 - 4 F_{0,3} F_{2,1} F_{3,0}^2}{2},\\
G_{0,4}&=
 F_{0,4} F_{3,0}^3 - 2 F_{0,3} F_{1,2} F_{3,0}^3.
\endaligned
\]

In $B_2$ when we assume $G_{0,3}\equiv0$, $G_{3,0}$ becomes a relative invariant $I^{\sf rel}_{3,0}$. In terms of the original Taylor coefficients
\begin{equation}\label{explicit-30}
\small
\aligned
I^{\sf rel}_{3,0}=&-\frac{1}{2 \sqrt{2} (F_{2,0} + 2 F_{1,1} + F_{0,2})^{\frac{3}{2}}(F_{1,1}^2 - F_{2,0} F_{0,2})^{\frac{3}{2}}} \\
&\times\Big\{-4 F_{0,3} F_{1,1}^3 + 6 F_{0,2} F_{1,1}^2 F_{1,2} + 3 F_{0,2} F_{0,3} F_{1,1} F_{2,0} - 
 6 F_{0,3} F_{1,1}^2 F_{2,0} - 3 F_{0,2}^2 F_{1,2} F_{2,0} + 9 F_{0,2} F_{1,1} F_{1,2} F_{2,0} \\
 &+ 
 3 F_{0,2} F_{0,3} F_{2,0}^2 - 3 F_{0,3} F_{1,1} F_{2,0}^2 + 9 F_{0,2} F_{1,2} F_{2,0}^2 + 
 3 F_{1,1} F_{1,2} F_{2,0}^2 - F_{0,3} F_{2,0}^3 + 4 F_{0,3} F_{1,1}^2 \sqrt{F_{1,1}^2 - F_{0,2} F_{2,0}} \\
 &- 
 6 F_{0,2} F_{1,1} F_{1,2} \sqrt{F_{1,1}^2 - F_{0,2} F_{2,0}} - 
 F_{0,2} F_{0,3} F_{2,0} \sqrt{F_{1,1}^2 - F_{0,2} F_{2,0}} + 
 6 F_{0,3} F_{1,1} F_{2,0} \sqrt{F_{1,1}^2 - F_{0,2} F_{2,0}} \\
 &- 
 9 F_{0,2} F_{1,2} F_{2,0} \sqrt{F_{1,1}^2 - F_{0,2} F_{2,0}} + 
 3 F_{0,3} F_{2,0}^2 \sqrt{F_{1,1}^2 - F_{0,2} F_{2,0}} + 
 3 F_{1,2} F_{2,0}^2 \sqrt{F_{1,1}^2 - F_{0,2} F_{2,0}} \\
 &- 3 F_{0,2}^2 F_{1,1} F_{2,1} - 
 9 F_{0,2}^2 F_{2,0} F_{2,1} - 9 F_{0,2} F_{1,1} F_{2,0} F_{2,1} - 6 F_{1,1}^2 F_{2,0} F_{2,1} + 
 3 F_{0,2} F_{2,0}^2 F_{2,1} \\
 &+ 3 F_{0,2}^2 \sqrt{F_{1,1}^2 - F_{0,2} F_{2,0}} F_{2,1} - 
 9 F_{0,2} F_{2,0} \sqrt{F_{1,1}^2 - F_{0,2} F_{2,0}} F_{2,1} - 
 6 F_{1,1} F_{2,0} \sqrt{F_{1,1}^2 - F_{0,2} F_{2,0}} F_{2,1} + F_{0,2}^3 F_{3,0} \\
 &+ 3 F_{0,2}^2 F_{1,1} F_{3,0} + 
 6 F_{0,2} F_{1,1}^2 F_{3,0} + 4 F_{1,1}^3 F_{3,0} - 3 F_{0,2}^2 F_{2,0} F_{3,0} - 
 3 F_{0,2} F_{1,1} F_{2,0} F_{3,0} \\
 &+ 3 F_{0,2}^2 \sqrt{F_{1,1}^2 - F_{0,2} F_{2,0}} F_{3,0} + 
 6 F_{0,2} F_{1,1} \sqrt{F_{1,1}^2 - F_{0,2} F_{2,0}} F_{3,0} + 
 4 F_{1,1}^2 \sqrt{F_{1,1}^2 - F_{0,2} F_{2,0}} F_{3,0} \\
 &- F_{0,2} F_{2,0} \sqrt{F_{1,1}^2 - F_{0,2} F_{2,0}} F_{3,0}\Big\}
\endaligned
\end{equation}
which is non-zero in this branch.

\subsection{Final loop for $B_2$}
The stabilizer of
\[
u=F(x,y)=xy+\frac{x^3}{6}+O(4)=xy+\frac{x^3}{6}+\sum\limits_{j+k\geqslant 4}\frac{F_{j,k}}{j!k!}x^jy^k,
\]
is
\[
\big\{x'=\mu x, y'=\mu^2 y, u'=\mu^3 u~|~\mu\in\bc^*\big\}.
\]
It acts on fourth order Taylor coefficients by sending $F_{4-k,k}$ to $\mu^{k+1}F_{4-k,k}$. So all fourth order Taylor coefficients $F_{4-k,k},k=0,1,2,3,4$ are relative invariants. We denote them by $I^{\sf rel}_{4-k,k},k=0,1,2,3,4$. They can be explicitly calculated by composing the two relations $G_{j,k}$ in terms of $F_{j,k}$ obtained in previous subsections. This gives
\[
\aligned
I^{\sf rel}_{4,0}&=\frac{P_{264}}{\sqrt{2}(F_{2,0}+2F_{1,1}+F_{0,2})^{\frac{1}{2}}(F_{1,1}^2 - F_{2,0}F_{0,2})^{\frac{3}{2}}P_{36}},\\
I^{\sf rel}_{3,1}&=\frac{P_{284}}{16(F_{2,0}+2F_{1,1}+F_{0,2})^2(F_{1,1}^2 - F_{2,0}F_{0,2})^3},\\
I^{\sf rel}_{2,2}&=-\frac{P_{26}P_{36}}{16\sqrt{2}(F_{2,0}+2F_{1,1}+F_{0,2})^{\frac{3}{2}}(F_{1,1}^2 - F_{2,0}F_{0,2})^{\frac{9}{2}}},\\
I^{\sf rel}_{1,3}&=\frac{Q_{284}P_{36}^2}{128(F_{2,0}+2F_{1,1}+F_{0,2})^5 (F_{1,1}^2 - F_{2,0}F_{0,2})^6},\\
I^{\sf rel}_{0,4}&=\frac{Q_{264}P_{36}^3}{64\sqrt{2}(F_{2,0}+2F_{1,1}+F_{0,2})^{\frac{13}{2}}(F_{1,1}^2 - F_{2,0}F_{0,2})^{\frac{15}{2}}},
\endaligned
\]
where $P_{36}$ is the numerator of $I^{\sf rel}_{3,0}$ and where $P_{264}$, $P_{284}$, $P_{26}$, $Q_{264}$, $Q_{284}$ are polynomials in $\bz[F_{j,k},\sqrt{F_{1,1}^2-F_{2,0}F_{0,2}}]$ having the indicated number of monomials. Moreover,
\[
P_{264}-Q_{264}=R_{156}, \ \ P_{284}-Q_{284}=4(F_{2,0}+2F_{1,1}+F_{0,2})\sqrt{F_{1,1}^2-F_{2,0}F_{0,2}} R_{66},
\]
where $R_{156}$ and $R_{66}$ are polynomials in $\bz[F_{j,k}]$ having the indicated number of monomials.

If $I^{\sf rel}_{4,0}\neq0$ we can normalize it to 1 by a unique choice of $\mu$. In this case we get the normal form
\[
u=xy+\frac{x^3}{6}+\frac{x^4}{24}+\sum\limits_{k=1}^4\frac{I_{4-k,k}}{(4-k)!k!}x^{4-k}y^k+O(5)
\]
where $I_{4-k,k}=\frac{I^{\sf rel}_{4-k,k}}{(I^{\sf rel}_{4,0})^{k+1}}$ are invariants.

Here we use Lie's principle: whenever we obtain a relative invariant, we only treat the cases where the relative invariant is non-zero or identically zero. This is true for generic points, i.e. points outside an analytic subset of codimension at least 1 on the concerned surface.

We may conclude our branching by a diagram
\[
\xymatrix{
\ar[rr]
\boxed{F_{1,1}^2-F_{2,0}F_{0,2}\neq0}
\ar[drr]
&&
{\sf Pick}\neq 0,
&&
&&
\\
&&
\ar[rr]
{\sf Pick}\equiv0
\ar[ddrr]
&&
\ar[rr]
I^{\sf rel}_{0,3}\equiv0\neq I^{\sf rel}_{3,0}
\ar[drr]
&&
I^{\sf rel}_{4,0}\neq0,
\\
&&
&&
&&
I^{\sf rel}_{4,0}\equiv0,
\\
&&
&&
I^{\sf rel}_{0,3}\equiv0\equiv I^{\sf rel}_{3,0}.
}
\]

The following sections study the existence of homogeneous models branch by branch.

\section{Non-vanishing Pick, branch $B_1$}\label{sect-nonvanishing-Pick}
According to \eqref{G-30-21-12-03}, every non-degenerate surface with non-vanishing Pick is $A(3)$-equivalent to a graph
\[
u=xy+\frac{x^3}{6}+\frac{y^3}{6}+\sum\limits_{j+k\geqslant 4}\frac{I_{j,k}}{j!k!}x^jy^k.
\]
The stabilizer group is discrete. By Fels-Olver's theory \cite[Thm 13.3]{Fels-Olver-1999}, all invariants are generated by the order 4 invariants $I_{j,4-j}$ for $j=0,1,2,3,4$ and their derivatives.

Olver's recurrence formulas, at order 4, are
\begin{equation}
\scriptsize
\aligned
\cd_xI_{4, 0} &= -8\,I_{1, 3}\,I_{4, 0}-160\,I_{4, 0}^2-144\,I_{3, 1}+I_{5, 0}, &\cd_yI_{4, 0} &= -32\,I_{0, 4}\,I_{4, 0}-40\,I_{3, 1}\,I_{4, 0}-48\,I_{2, 2}+I_{4, 1}+216, \\
 \cd_xI_{3, 1} &= -4\,I_{1, 3}\,I_{3, 1}-32\,I_{3, 1}\,I_{4, 0}-84\,I_{2, 2}+I_{4, 1}+216, &\cd_yI_{3, 1} &= -16\,I_{0, 4}\,I_{3, 1}-8\,I_{3, 1}^2-72\,I_{1, 3}-72\,I_{4, 0}+I_{3, 2}, \\
  \cd_xI_{2, 2} &= -4\,I_{1, 3}\,I_{2, 2}-16\,I_{2, 2}\,I_{4, 0}-36\,I_{1, 3}+I_{3, 2}, &\cd_yI_{2, 2} &= -16\,I_{0, 4}\,I_{2, 2}-4\,I_{2, 2}\,I_{3, 1}-36\,I_{3, 1}+I_{2, 3}, \\
   \cd_xI_{1, 3} &= -8\,I_{1, 3}^2-16\,I_{1, 3}\,I_{4, 0}-72\,I_{0, 4}-72\,I_{3, 1}+I_{2, 3}, &\cd_yI_{1, 3} &= -32\,I_{0, 4}\,I_{1, 3}-4\,I_{1, 3}\,I_{3, 1}-84\,I_{2, 2}+I_{1, 4}+216, \\
    \cd_xI_{0, 4} &= -40\,I_{0, 4}\,I_{1, 3}-32\,I_{0, 4}\,I_{4, 0}-48\,I_{2, 2}+I_{1, 4}+216, &\cd_yI_{0, 4} &= -160\,I_{0, 4}^2-8\,I_{0, 4}\,I_{3, 1}-144\,I_{1, 3}+I_{0, 5}.
\endaligned
\end{equation}
There are 10 equations. One can solve the 6 order 5 invariants $I_{j,5-j}$ for $j=0,1,2,3,4,5$ in terms of the order 4 invariants and their derivatives.

For simplicity, we say two functions $R(I_{j_1,k_1},\cd_x^{j_2}\cd_y^{k_2}I_{j_3,k_3})$ and $\tilde{R}$ are {\sl conjugate} if
\[
\tilde{R}(I_{k_1,j_1},\cd_x^{k_2}\cd_y^{j_2}I_{k_3,j_3})=R(I_{j_1,k_1},\cd_x^{j_2}\cd_y^{k_2}I_{j_3,k_3}),
\]
i.e. after switching $x$ and $y$, they are the same. Two equations are {\sl conjugate} if they can be written as $0=R$ and $0=\tilde{R}$ for some conjugate pair $R$ and $\tilde{R}$. For example, the 10 recurrence formulas of order 4 contain 5 conjugate pairs of equations.

From
\[
\aligned
\cd_yI_{4, 0} &= -32\,I_{0, 4}\,I_{4, 0}-40\,I_{3, 1}\,I_{4, 0}-48\,I_{2, 2}+I_{4, 1}+216,\\
\cd_xI_{3, 1} &= -4\,I_{1, 3}\,I_{3, 1}-32\,I_{3, 1}\,I_{4, 0}-84\,I_{2, 2}+I_{4, 1}+216,
\endaligned
\]
we can eliminate $I_{4,1}$ and solve $I_{2,2}$
\[
I_{2, 2} = \tfrac{8}{9}\,I_{4, 0}\,I_{0, 4}-\tfrac{1}{9}\,I_{1, 3}\,I_{3, 1}+\tfrac{2}{9}\,I_{4, 0}\,I_{3, 1}-\tfrac{1}{36}\,\cd_xI_{3, 1}+\tfrac{1}{36}\,\cd_yI_{4, 0}.
\]
From the two conjugate equations we get a conjugate solution
\[
I_{2, 2} = \tfrac{2}{9}\,I_{1, 3}\,I_{0, 4}+\tfrac{8}{9}\,I_{4, 0}\,I_{0, 4}-\tfrac{1}{9}\,I_{1, 3}\,I_{3, 1}+\tfrac{1}{36}\,\cd_xI_{0, 4}-\tfrac{1}{36}\,\cd_yI_{1, 3}.
\]
To conclude, $I_{2,2}$ can be solved in terms of the other 4 invariants and their derivatives.

Under extra assumptions on genericity, for example
\[
\det\left(\!
\begin{array}{cc}
\cd_x I_{j,4-j} & \cd_y I_{j,4-j}\\
\cd_x^2 I_{j,4-j} & \cd_y \,\cd_x I_{j,4-j}
\end{array}
\!\right)\neq0,
\]
for some $0\leqslant j\leqslant 4$, one may find generating systems with fewer elements by investigating the Lie bracket
\[
[\cd_x,\cd_y]=(\tfrac{2}{3}\,I_{3, 1}+\tfrac{4}{3}\,I_{0, 4})\,\cd_x+(-\tfrac{4}{3}\,I_{4, 0}-\tfrac{2}{3}\,I_{1, 3})\,\cd_y
\]
and by using the same method in \cite{Olver-2007, Arnaldsson-Valiquette-2020}. But it is not the case for homogeneous surfaces, where $\cd_x I_{j,4-j}=\cd_y I_{j,4-j}=0$ since $I_{j,4-j}$ are constant.

For homogeneous surfaces, all invariants have to be constant. Thus all left hand sides in these recurrence formulas are 0. The over-determined linear system is solvable if and only if 4 more conditions among $I_{j,4-j}$ are satisfied:
\[
\aligned
(E1) \ \ 0&=8\,I_{0, 4}\,I_{4, 0}-I_{1, 3}\,I_{3, 1}+2\,I_{3, 1}\,I_{4, 0}-9\,I_{2, 2},\\
(E2) \ \ 0&= 2\,I_{0, 4}\,I_{1, 3}+8\,I_{0, 4}\,I_{4, 0}-I_{1, 3}\,I_{3, 1}-9\,I_{2, 2},\\
(E3) \ \ 0&= 4\,I_{0, 4}\,I_{3, 1}-I_{1, 3}\,I_{2, 2}-4\,I_{2, 2}\,I_{4, 0}+2\,I_{3, 1}^2+9\,I_{1, 3}+18\,I_{4, 0},\\
(E4) \ \ 0&= 4\,I_{0, 4}\,I_{2, 2}-2\,I_{1, 3}^2-4\,I_{1, 3}\,I_{4, 0}+I_{2, 2}\,I_{3, 1}-18\,I_{0, 4}-9\,I_{3, 1}.
\endaligned
\]
In the equation $(E1)$ we solve
\begin{equation}\label{sol22o}
I_{2,2}=\frac{8}{9}I_{0,4}I_{4,0}-\frac{1}{9}I_{1,3}I_{3,1}+\frac{2}{9}I_{3,1}I_{4,0}.
\end{equation}
We put this solution back to $(E2)$, $(E3)$, $(E4)$
\[
\aligned
(F1) \ \ 0&=I_{1,3}\,I_{0,4}-I_{4,0}\,I_{3,1},\\
(F2) \ \ 0&=4\,I_{0, 4}\,I_{3, 1}-\tfrac{8}{9}\,I_{0, 4}\,I_{1, 3}\,I_{4, 0}+\tfrac{1}{9}\,I_{1, 3}^2\,I_{3, 1}+\tfrac{2}{9}\,I_{1, 3}\,I_{3, 1}\,I_{4, 0}\\
& \ \ \ -\tfrac{32}{9}\,I_{0, 4}\,I_{4, 0}^2-\tfrac{8}{9}\,I_{3, 1}\,I_{4, 0}^2+2\,I_{3, 1}^2+9\,I_{1, 3}+18\,I_{4, 0},\\
(F3) \ \ 0&=\tfrac{32}{9}\,I_{0, 4}^2\,I_{4, 0}-\tfrac{4}{9}\,I_{0, 4}\,I_{1, 3}\,I_{3, 1}+\tfrac{16}{9}\,I_{0, 4}\,I_{3, 1}\,I_{4, 0}-2\,I_{1, 3}^2\\
& \ \ \ -4\,I_{1, 3}\,I_{4, 0}-\tfrac{1}{9}\,I_{1, 3}\,I_{3, 1}^2+\tfrac{2}{9}\,I_{3, 1}^2\,I_{4, 0}-18\,I_{0, 4}-9\,I_{3, 1}.
\endaligned
\]
\noindent {\bf Case 1:} If $I_{4,0}\neq0$ we solve $I_{3,1}=\frac{I_{0,4}I_{1,3}}{I_{4,0}}$. Replace $I_{2,2}$ and $I_{3,1}$ in $(E3)$ and $(E4)$ we get
\[
\aligned
(G1) \ \ 0&=\frac{2I_{4,0}+I_{1,3}}{9I_{4,0}^2}\Big(I_{0,4}I_{1,3}^2I_{4,0}-16I_{0,4}I_{4,0}^3+18I_{0,4}^2I_{1,3}+81I_{4,0}^2\Big),\\
(G2) \ \ 0&=\frac{2I_{4,0}+I_{1,3}}{9I_{4,0}^2}\Big(I_{0,4}^2I_{1,3}^2-16I_{0,4}^2I_{4,0}^2+18I_{1,3}I_{4,0}^2+81I_{4,0}I_{0,4}\Big).
\endaligned
\]
\noindent {\bf Subcase 1-1:} If $I_{1,3}=-2I_{4,0}$ then all relations are satisfied. If we write $a:=I_{4,0}$, $b:=I_{0,4}$ then
\[
I_{3,1}=-2b, \ \ I_{1,3}=-2a, \ \ I_{2,2}=0, \ \ a\in\bc^*, b\in\bc.
\]
The homogeneous candidate is
\[
u=xy+\frac{x^3}{6}+\frac{y^3}{6}+\frac{a}{24}x^4-\frac{b}{3}x y^3-\frac{a}{3}x^3y+\frac{b}{24}y^4+O(5), \ \ a\in\bc^*,b\in\bc,
\]
and it corresponds to N1 in \cite{Eastwood-Ezhov-1999} with $a\in\bc^*$, $b\in\bc$.

\noindent {\bf Subcase 1-2:} $I_{1,3}\neq-2I_{4,0}$. Then $(G1)$ and $(G2)$ become
\[
\aligned
(H1) \ \ 0&=I_{0,4}I_{1,3}^2I_{4,0}-16I_{0,4}I_{4,0}^3+18I_{0,4}^2I_{1,3}+81I_{4,0}^2,\\
(H2) \ \ 0&=I_{0,4}^2I_{1,3}^2-16I_{0,4}^2I_{4,0}^2+18I_{1,3}I_{4,0}^2+81I_{4,0}I_{0,4}.
\endaligned
\]
We calculate the two sides of $I_{0,4}\,(H1)-I_{4,0}\,(H2)$
\[
0=18\,I_{1,3}\,(I_{0,4}^3-I_{4,0}^3).
\]

\noindent{\bf Subsubcase 1-2-1:} If $I_{1,3}=0$ then by $(F1)$, $I_{3,1}=0$. The equations $(H1)$, $(H2)$ become
\[
\aligned
0&=-16I_{0,4}I_{4,0}^3+81I_{4,0}^2,\\
0&=-16I_{0,4}^2I_{4,0}^2+81I_{4,0}I_{0,4},
\endaligned
\]
and we solve $I_{0,4}=\frac{81}{16\,I_{4,0}}$. Put back to \eqref{sol22o} we get $I_{2,2}=\tfrac{9}{2}$. That is N4 in \cite{Eastwood-Ezhov-1999}.

\noindent{\bf Subsubcase 1-2-2:} If $I_{1,3}\neq0$ then $I_{0,4}=I_{4,0}\,\omega^j$ for some $j=0,1,2$. Here we recall $\omega=e^{2\pi i/3}$ the cube root of unity. By a transformation in $\mathcal{G}_0$, they are equivalent to $I_{0,4}=I_{4,0}\in\bc^*$ which implies $I_{3,1}=I_{1,3}$ by $(F1)$. Put it back to $(H1)$, $(H2)$ we get the same equation
\[
\aligned
0&=I_{1,3}^2\,I_{4,0}^2-16\,I_{4,0}^4+18\,I_{1,3}\,I_{4,0}^2+81\,I_{4,0}^2\\
&=I_{4,0}^2\,(I_{1,3}+9-4\,I_{4,0})\,(I_{1,3}+9\,I_{4,0}+4\,I_{4,0}).
\endaligned
\]
Thus either
\[
I_{3,1}=I_{1,3}=4\,I_{4,0}-9, \ \ I_{2,2}=6\,I_{4,0}-9, \ \ I_{0,4}=I_{4,0}\in\bc^*,
\]
corresponds to N3 in \cite{Eastwood-Ezhov-1999} with $b\in\bc^*$, or
\[
I_{3,1}=I_{1,3}=-4\,I_{4,0}-9, \ \ I_{2,2}=-\tfrac{16}{9}\,I_{4,0}^2-10\,I_{4,0}-9, \ \ I_{0,4}=I_{4,0}\in\bc^*,
\]
corresponds to N2 in \cite{Eastwood-Ezhov-1999} with $b\in\bc^*$.

\noindent{\bf Case 2:} If $I_{4,0}=0$, then $I_{1,3}I_{0,4}=0$ from $(F1)$.

\noindent{\bf Subcase 2-1:} If $I_{1,3}=0$, then $(F2)$ and $(F3)$ becomes
\[
0 = (2 I_{0,4} + I_{3,1} )  I_{3,1}, \ \ 0 =  2 I_{0,4} + I_{3,1}.
\]
The solution is $I_{0,4}=b\in\bc$, $I_{3,1}=-2b$, corresponds to N1 in \cite{Eastwood-Ezhov-1999} with $a=0$, $b\in\bc$.

\noindent{\bf Subcase 2-2:} If $I_{1,3}\neq0$ then $I_{0,4}=0$, and $(F2)$, $(F3)$ become
\[
\aligned
(H1') \ \ 0&=2 I_{3,1}^2+\tfrac{1}{9}I_{1,3}^2I_{3,1}+9I_{1,3},\\
(H2') \ \ 0&=-\tfrac{1}{9}I_{1,3}I_{3,1}^2-9I_{3,1}-2I_{1,3}^2.
\endaligned
\]
We calculate the two sides of $I_{3,1}\,(H1')+I_{1,3}\,(H2')$
\[
0=2\,(I_{3,1}^3-I_{1,3}^3)
\]
So $I_{3,1}=I_{1,3}\,\omega^j$ for some $j=0,1,2$. By a transformation in $\mathcal{G}_0$, they are equivalent to $I_{3,1}=I_{1,3}=t\in\bc^*$. Put it back to $(H1')$, $(H2')$ we get the same equation
\[
\aligned
0&=2t^2+\tfrac{1}{9}t^3+9t, \ \ t\in\bc^*\\
&=\frac{t}{9}(t+9)^2.
\endaligned
\]
So $I_{3,1}=I_{1,3}=t=-9$. By \eqref{sol22o}, $I_{2,2}=-9$. The solution
\[
I_{3,1}=I_{1,3}=I_{2,2}=-9, \ \ I_{4,0}=I_{0,4}=0,
\]
corresponds to N2 and N3 in \cite{Eastwood-Ezhov-1999} with $b=0$.

\section{Vanishing $I^{\sf rel}_{0,3}$ but $I^{\sf rel}_{4,0}\neq0$, branch $B_{2\cdot1}$}\label{sect-v03-nv40}
Suppose $I^{\sf rel}_{0,3}=0$, by homogeneity $I^{\sf rel}_{0,3}\equiv0$. Thus we may solve, on the jet space $J_{2,1}^3$
\[
\aligned
u_{0,3}=&\frac{1}{(u_{1,1}+u_{2,0}+\sqrt{u_{1,1}^2-u_{2,0}u_{0,2}})^3}\\
&\times\Big\{6 u_{0,2} u_{1,1}^2 u_{1,2} - 3 u_{0,2}^2 u_{1,2} u_{2,0} + 9 u_{0,2} u_{1,1} u_{1,2} u_{2,0} + 
 9 u_{0,2} u_{1,2} u_{2,0}^2 + 3 u_{1,1} u_{1,2} u_{2,0}^2 \\
 &+ 
 6 u_{0,2} u_{1,1} u_{1,2} \sqrt{u_{1,1}^2 - u_{0,2} u_{2,0}} + 
 9 u_{0,2} u_{1,2} u_{2,0} \sqrt{u_{1,1}^2 - u_{0,2} u_{2,0}} - 
 3 u_{1,2} u_{2,0}^2 \sqrt{u_{1,1}^2 - u_{0,2} u_{2,0}} \\
 &- 3 u_{0,2}^2 u_{1,1} u_{2,1} - 
 9 u_{0,2}^2 u_{2,0} u_{2,1} - 9 u_{0,2} u_{1,1} u_{2,0} u_{2,1} - 6 u_{1,1}^2 u_{2,0} u_{2,1} + 
 3 u_{0,2} u_{2,0}^2 u_{2,1} \\
 &- 3 u_{0,2}^2 \sqrt{u_{1,1}^2 - u_{0,2} u_{2,0}} u_{2,1} + 
 9 u_{0,2} u_{2,0} \sqrt{u_{1,1}^2 - u_{0,2} u_{2,0}} u_{2,1} + 
 6 u_{1,1} u_{2,0} \sqrt{u_{1,1}^2 - u_{0,2} u_{2,0}} u_{2,1} \\
 &+ u_{0,2}^3 u_{3,0} + 3 u_{0,2}^2 u_{1,1} u_{3,0} + 
 6 u_{0,2} u_{1,1}^2 u_{3,0} + 4 u_{1,1}^3 u_{3,0} - 3 u_{0,2}^2 u_{2,0} u_{3,0} - 
 3 u_{0,2} u_{1,1} u_{2,0} u_{3,0} \\
 &- 3 u_{0,2}^2 \sqrt{u_{1,1}^2 - u_{0,2} u_{2,0}} u_{3,0} - 
 6 u_{0,2} u_{1,1} \sqrt{u_{1,1}^2 - u_{0,2} u_{2,0}} u_{3,0} - 
 4 u_{1,1}^2 \sqrt{u_{1,1}^2 - u_{0,2} u_{2,0}} u_{3,0} \\
 &+ u_{0,2} u_{2,0} \sqrt{u_{1,1}^2 - u_{0,2} u_{2,0}} u_{3,0}\Big\},
\endaligned
\]
By taking total differentials $D_x$ and $D_y$, we may solve $u_{0,4}$ and $u_{1,3}$ on $J_{2,1}^4$. Indeed we can solve $u_{n-3-k,3+k}$ for any $0\leqslant k\leqslant n-3$ on $J_{2,1}^n$ in terms of $u_{j,0}$, $u_{j-1,1}$ and $u_{j-2,2}$ for $0\leqslant j\leqslant n$. We call those $u_{n-3-k,3+k}$ as {\sl dependent jet coordinates} and the others as {\sl independent jet coordinates}.

\begin{Def} The {\sl subjet} of $I^{\sf rel}_{0,3}\equiv0$ is a series of submanifolds $SJ_{2,1}^n\subset J_{2,1}^n$ for $n\geqslant 3$ determined by the PDEs $D_x^jD_y^k (I^{\sf rel}_{0,3})=0$ with $0\leqslant j+k\leqslant n-3$.
\end{Def}

\begin{Prop} The germ of an analytic surface lies in $SJ_{2,1}^n$ for any $n\geqslant 3$ if and only if the surface has $I^{\sf rel}_{0,3}\equiv0$.
\end{Prop}

\begin{Prop} In a neighborhood of the cross-section corresponding to the normalization
\[
\aligned
&u_{0,3}=0, \\
&u_{0,2}=0, \ \ u_{1,2}=0, \\
&u_{0,1}=0, \ \ u_{1,1}=1, \ \ u_{2,1}=0, \\
&u_{0,0}=0, \ \ u_{1,0}=0, \ \ u_{2,0}=0, \ \ u_{3,0}=1,\endaligned
\]
the subjet $SJ_{2,1}^n$ can be graphed by the solutions $u_{m-3-k,3+k}=\mathcal{R}_{m-3-k,3+k}(u_{j,0},u_{j-1,1},u_{j-2,2})$ for $0\leqslant k\leqslant m-3$ and $3\leqslant m\leqslant n$ obtained above.
\end{Prop}

There is a natural projection $\pi^n$ from $SJ_{2,1}^n$ to the span of independent coordinates. The prolonged group action and its infinitesimal generators on $J_{2,1}^n$, restricted to the invariant submanifold $SJ_{2,1}^n$ (invariant because $I^{\sf rel}_{3,0}=0$ is an invariant property), can be pushed forward to the space of independent coordinates. Thus we can run Olver's recurrence formulas on the span and get relations among invariants associated to independent coordinates.

In this branch $B_{2\cdot1}$ we assume $I^{\sf rel}_{4,0}\neq0$, so we normalize it to 1. The normal form is
\[
u=xy+\frac{x^3}{6}+\frac{x^4}{24}+\frac{I_{3,1}}{6}x^3y+\frac{I_{2,2}}{4}x^2y^2+\frac{I_{1,3}}{6}xy^3+\frac{I_{0,4}}{24}y^4+O(5).
\]
There are 2 invariants $I_{3,1}$ and $I_{2,2}$ of order 4 from the independent coordinates. The other 2 invariants $I_{1,3}$, $I_{0,4}$ are dependent because if we solve $I^{\sf rel}_{0,3}=0$ for this power series, we get $I_{1,3}=0$ and $I_{0,4}=0$.

There are 3 order 5 invariants $I_{5,0}$, $I_{4,1}$ and $I_{3,2}$. Again by Fels-Olver's theory \cite[Thm 13.3]{Fels-Olver-1999}, all invariants are generated by $I_{3,1}$, $I_{2,2}$, $I_{5,0}$, $I_{4,1}$ and their derivatives.

Olver's recurrence formulas, at order 4, are
\[
\aligned
\cd_x I_{3,1}& =I_{4,1}+8I_{3,1}^2-\tfrac{7}{2}I_{2,2}+2I_{3,1}-2I_{5,0}I_{3,1},\\
\cd_y I_{3,1}&=4I_{3,1}I_{2,2}+2I_{3,1}^2-2I_{4,1}I_{3,1}+I_{3,2},\\
\cd_x I_{2,2}&=12I_{3,1}I_{2,2}-3I_{5,0}I_{2,2}+4I_{2,2}+I_{3,2},\\
\cd_y I_{2,2}&=6I_{2,2}^2+4I_{3,1}I_{2,2}-3I_{4,1}I_{2,2}.
\endaligned
\]
In the first formula, we can solve
\[
I_{4, 1} = -8\,I_{3, 1}^2+2\,I_{5, 0}\,I_{3, 1}+\cd_x\,I_{3, 1}+\tfrac{7}{2}\,I_{2, 2}-2\,I_{3, 1}.
\]
Thus $I_{3,1}$, $I_{2,2}$, $I_{5,0}$ are generators.

Like in the branch $B_1$, under extra assupmtions on genericity, for example
\[
\det\left(\!
\begin{array}{cc}
\cd_x I_{3,1} & \cd_y I_{3,1}\\
\cd_x^2 I_{3,1} & \cd_y \cd_x I_{3,1}
\end{array}
\!\right)\neq0,
\]
one may find generating systems with fewer elements by investigating the Lie bracket
\[
[\cd_x,\cd_y]=(-I_{3, 1}-2\,I_{2,2}+I_{4,1})\,\cd_x+(8\,I_{3,1}+3-2\,I_{5,0})\,\cd_y
\]
and by using the same method in \cite{Olver-2007, Arnaldsson-Valiquette-2020}. But again it is not the case for homogeneous surfaces, where $\cd_x I_{j,4-j}=\cd_y I_{j,4-j}=0$.

Homogeneous surfaces always have constant invariants, i.e. $\cd_x I_{j,k}=\cd_y I_{j,k}=0$. So for homogeneous surfaces the formulas are
\[
\aligned
(E41) \ \ 0& =I_{4,1}+8I_{3,1}^2-\tfrac{7}{2}I_{2,2}+2I_{3,1}-2I_{5,0}I_{3,1},\\
(E42) \ \ 0&=4I_{3,1}I_{2,2}+2I_{3,1}^2-2I_{4,1}I_{3,1}+I_{3,2},\\
(E43) \ \ 0&=12I_{3,1}I_{2,2}-3I_{5,0}I_{2,2}+4I_{2,2}+I_{3,2},\\
(E44) \ \ 0&=6I_{2,2}^2+4I_{3,1}I_{2,2}-3I_{4,1}I_{2,2}.
\endaligned
\]
We may solve from $(E41)$, $(E42)$
\begin{equation}\label{sol5}
\aligned
I_{4,1} &= -8\,I_{3, 1}^2+2\,I_{5, 0}\,I_{3, 1}+\tfrac{7}{2}\,I_{2, 2}-2\,I_{3, 1},  \\
I_{3,2} &= -16\,I_{3, 1}^3+4\,I_{3, 1}^2\,I_{5, 0}+3\,I_{2, 2}\,I_{3, 1}-6\,I_{3, 1}^2.
\endaligned
\end{equation}
We put them back to $(E43)$, $(E44)$
\[
\aligned
(F41) \ \ 0&=-16\,I_{3, 1}^3+4\,I_{3, 1}^2\,I_{5, 0}+15\,I_{2, 2}\,I_{3, 1}-3\,I_{2, 2}\,I_{5, 0}-6\,I_{3, 1}^2+4\,I_{2, 2}, \\
(F42) \ \ 0&=-\tfrac{I_{2,2}}{2}\,(-48\,I_{3, 1}^2+12\,I_{3, 1}\,I_{5, 0}+9\,I_{2, 2}-20\,I_{3, 1}).
\endaligned
\]

The recurrence formulas of order 5 provide
\[
\aligned
\cd_x I_{5,0}&=I_{6, 0}+8I_{5, 0}I_{3, 1}-\tfrac{5}{2}I_{4, 1}+I_{5, 0}-\tfrac{5}{2}I_{3, 1}-5I_{2, 2}-2I_{5, 0}^2, \\
\cd_y I_{5,0}&=I_{5, 1}+I_{5, 0}I_{3, 1}+4I_{5, 0}I_{2, 2}-\tfrac{5}{2}I_{2, 2}-2I_{5, 0}I_{4, 1},\\
\cd_xI_{4, 1} &= I_{5, 1}+12I_{3, 1}I_{4, 1}-2I_{3, 2}+3I_{4, 1}-4I_{3, 1}^2-\tfrac{5}{2}I_{2, 2}-3I_{5, 0}I_{4, 1},\\
\cd_yI_{4, 1} &= -4I_{2, 2}I_{3, 1}+6I_{2, 2}I_{4, 1}+3I_{3, 1}I_{4, 1}-3I_{4, 1}^2+I_{4, 2},\\
\cd_xI_{3, 2} &= I_{4, 2}+16I_{3, 1}I_{3, 2}+5I_{3, 2}-\tfrac{17}{2}I_{3, 1}I_{2, 2}-4I_{5, 0}I_{3, 2},\\
\cd_yI_{3, 2} &= 8I_{2, 2}I_{3, 2}+5I_{3, 1}I_{3, 2}-4I_{3, 2}I_{4, 1}.\\
\endaligned
\]
In the homogeneous case the left hand sides are all 0. Replacing $I_{4,1}$, $I_{3,2}$ by their solutions above, the 6 equations become
\[
\aligned
(E51) \ \ 0&=I_{6, 0}+8I_{5, 0}I_{3, 1}-\tfrac{5}{2}I_{4, 1}+I_{5, 0}-\tfrac{5}{2}I_{3, 1}-5I_{2, 2}-2I_{5, 0}^2, \\
(E52) \ \ 0&=I_{5, 1}+I_{5, 0}I_{3, 1}+4I_{5, 0}I_{2, 2}-\tfrac{5}{2}I_{2, 2}-2I_{5, 0}I_{4, 1},\\
(E53) \ \ 0&= I_{5, 1}+12I_{3, 1}I_{4, 1}-2I_{3, 2}+3I_{4, 1}-4I_{3, 1}^2-\tfrac{5}{2}I_{2, 2}-3I_{5, 0}I_{4, 1},\\
(E54) \ \ 0&= -4I_{2, 2}I_{3, 1}+6I_{2, 2}I_{4, 1}+3I_{3, 1}I_{4, 1}-3I_{4, 1}^2+I_{4, 2},\\
(E55) \ \ 0&= I_{4, 2}+16I_{3, 1}I_{3, 2}+5I_{3, 2}-\tfrac{17}{2}I_{3, 1}I_{2, 2}-4I_{5, 0}I_{3, 2},\\
(E56) \ \ 0&= 8I_{2, 2}I_{3, 2}+5I_{3, 1}I_{3, 2}-4I_{3, 2}I_{4, 1}.\\
\endaligned
\]
We solve from $(E51)$, $(E52)$, $(E55)$ while using the solutions \eqref{sol5}
\begin{equation}\label{sol6}
\aligned
I_{6, 0} &= -20\,I_{3, 1}^2-3\,I_{5, 0}\,I_{3, 1}+2\,I_{5, 0}^2+\tfrac{55}{4}\,I_{2, 2}-\tfrac{5}{2}\,I_{3, 1}-I_{5, 0},\\
I_{5, 1} &= -16\,I_{3, 1}^2\,I_{5, 0}+4\,I_{5, 0}^2\,I_{3, 1}+3\,I_{2, 2}\,I_{5, 0}-5\,I_{5, 0}\,I_{3, 1}+\tfrac{5}{2}\,I_{2, 2},\\
I_{4, 2} &= 256\,I_{3, 1}^4-128\,I_{3, 1}^3\,I_{5, 0}+16\,I_{5, 0}^2\,I_{3, 1}^2-48\,I_{2, 2}\,I_{3, 1}^2\\
& \ \ \ +12\,I_{2, 2}\,I_{5, 0}\,I_{3, 1}+176\,I_{3, 1}^3-44\,I_{3, 1}^2\,I_{5, 0}-\tfrac{13}{2}\,I_{2, 2}\,I_{3, 1}+30\,I_{3, 1}^2.
\endaligned
\end{equation}
We put the solutions \eqref{sol5} and \eqref{sol6} back into $(E53)$, $(E54)$, $(E56)$
\[
\aligned
(F51) \ \ 0&=24\,I_{3, 1}^2\,I_{5, 0}-2\,I_{5, 0}^2\,I_{3, 1}-\tfrac{15}{2}\,I_{2, 2}\,I_{5, 0}+7\,I_{5, 0}\,I_{3, 1}\\
& \ \ \ +\tfrac{21}{2}\,I_{2, 2}-64\,I_{3, 1}^3+36\,I_{2, 2}\,I_{3, 1}-40\,I_{3, 1}^2-6\,I_{3, 1},\\
(F52) \ \ 0&=30\,I_{2, 2}\,I_{3, 1}+72\,I_{2, 2}\,I_{3, 1}^2-18\,I_{2, 2}\,I_{5, 0}\,I_{3, 1}-\tfrac{63}{4}\,I_{2, 2}^2\\
& \ \ \ +56\,I_{3, 1}^3-14\,I_{3, 1}^2\,I_{5, 0}+12\,I_{3, 1}^2+64\,I_{3, 1}^4-32\,I_{3, 1}^3\,I_{5, 0}+4\,I_{5, 0}^2\,I_{3, 1}^2,\\
(F53) \ \ 0&=-I_{3, 1}\,(-16\,I_{3, 1}^2+4\,I_{3, 1}\,I_{5, 0}+3\,I_{2, 2}-6\,I_{3, 1})\,(-32\,I_{3, 1}^2+8\,I_{3, 1}\,I_{5, 0}+6\,I_{2, 2}-13\,I_{3, 1}).
\endaligned
\]
We obtain necessary conditions for being homogeneous: $(F41)$, $(F42)$, $(F51)$, $(F52)$, $(F53)$.

\medskip
\noindent {\bf Case 1:} If $I_{2,2}=0$, then $(F42)$ is satisfied. The other equations become
\[
\aligned
(G11) \ \ 0&= -2\,I_{3, 1}^2\,(8\,I_{3, 1}-2\,I_{5, 0}+3),\\
(G12) \ \ 0&= -I_{3, 1}\,(8\,I_{3, 1}-I_{5, 0}+2)\,(8\,I_{3, 1}-2\,I_{5, 0}+3),\\
(G13) \ \ 0&= 2\,I_{3, 1}^2\,(4\,I_{3, 1}-I_{5, 0}+2)\,(8\,I_{3, 1}-2\,I_{5, 0}+3),\\
(G14) \ \ 0&= -2\,I_{3, 1}^3\,(8\,I_{3, 1}-2\,I_{5, 0}+3)\,(32\,I_{3, 1}-8\,I_{5, 0}+13).
\endaligned
\]
Thus either $I_{3,1}=0$, which corresponds to N6 in \cite{Eastwood-Ezhov-1999}, or $8\,I_{3, 1}-2\,I_{5, 0}+3=0$, which corresponds to N5 in \cite{Eastwood-Ezhov-1999}.

\medskip
\noindent {\bf Case 2:} If $I_{2,2}\neq0$, then $(F42)$ becomes
\[
0=-48\,I_{3, 1}^2+12\,I_{3, 1}\,I_{5, 0}+9\,I_{2, 2}-20\,I_{3, 1},
\]
where we can solve
\begin{equation}\label{sol22}
I_{2, 2} = \tfrac{16}{3}\,I_{3, 1}^2-\tfrac{4}{3}\,I_{5, 0}\,I_{3, 1}+\tfrac{20}{9}\,I_{3, 1}.
\end{equation}
We put this solution \eqref{sol22} back to $(F53)$
\[
0=-\tfrac{9}{2}\,I_{3,1}^3.
\]
Thus $I_{3,1}=0$. Put it back to \eqref{sol22} we get $I_{2,2}=0$ and we return to Case 1.

\section{Vanishing $I^{\sf rel}_{0,3}\equiv0\equiv I^{\sf rel}_{4,0}$, branch $B_{2\cdot2}$}\label{sect-v03-v40}

By $I^{\sf rel}_{0,3}\equiv0$ we can solve all $u_{\geqslant0,\geqslant3}$. By $I^{\sf rel}_{4,0}\equiv0$ we can solve all $u_{\geqslant4,\geqslant0}$. Thus only finitely many jet coordinates are independent, namely $u_{\leqslant 3,\leqslant 2}$. In the previous section we have already normalized 
\[
\aligned
u_{0,2}&=0,  & u_{1,2}&=0,  & u_{2,2}& ,& u_{3,2}& ,\\
u_{0,1}&=0,  & u_{1,1}&=1,  & u_{2,1}&=0 , &u_{3,1}& ,\\
u_{0,0}&=0,  & u_{1,1}&=0,  & u_{2,0}&=0 , &u_{3,0}&=1,
\endaligned
\]
only $u_{2,2}$, $u_{3,1}$, $u_{3,2}$ remain free.

However, the infinite PDE system of $D_x^j D_y^k I^{\sf rel}_{0,3}=0$ and $D_x^j D_y^k I^{\sf rel}_{4,0}=0$ is not always compatible. The compatibility condition is necessary for a surface to be homogeneous.

From $D_x^j D_y^k I^{\sf rel}_{0,3}=0$ for $j+k\leqslant 4$ we solve all $u_{j,k+3}$. At the cross-section
\[
u_{1,1}=u_{3,0}=1,  \ \ u_{0,0}=u_{1,0}=u_{0,1}=u_{2,0}=u_{0,2}=u_{2,1}=u_{1,2}=0,
\]
we have
\[
u_{3,3}=\frac{9}{2}\,u_{2,2}^2, \ \ u_{4,3}=15\,u_{3,2}\,u_{2,2},
\]
and all the other $u_{j,k+3}=0$. 

From $D_x^{j'} D_y^{k'} I^{\sf rel}_{4,0}=0$ for $j'+k'\leqslant 3$ we solve all $u_{j'+4,k'}$. At the same cross-section we have
\[
\aligned
u_{4,3}=&-6\,u_{0, 3}\,u_{3, 1}^2+6\,u_{0, 3}\,u_{1, 3}-6\,u_{0, 3}\,u_{3, 2}-6\,u_{0, 4}\,u_{3, 1}-12\,u_{1, 3}\,u_{2, 2}\\
&+12\,u_{2, 2}\,u_{3, 2}+8\,u_{2, 3}\,u_{3, 1}-u_{0, 5}+2\,u_{2, 4}.
\endaligned
\]
If the system is compatible, we may replace $u_{j,k+3}$ from both sides by solutions of $D_x^j D_y^k I^{\sf rel}_{0,3}=0$:
\[
15\,u_{3,2}\,u_{2,2}=12\,u_{3,2}\,u_{2,2}.
\]
Thus the PDE system of $D_x^j D_y^k I^{\sf rel}_{0,3}=0$ with $j+k\leqslant 4$ and $D_x^{j'} D_y^{k'} I^{\sf rel}_{4,0}=0$ with $j'+k'\leqslant 3$ is compatible if and only if $u_{3,2}\,u_{2,2}=0$. Either $u_{3,2}=0$ or $u_{2,2}=0$.

Furthermore, by checking $u_{5,3}$ we get
\[
u_{5,3}=75\,u_{2,2}^2\,u_{3,1}+12\,u_{3,2}^2=\frac{135}{2}\,u_{2,2}^2\,u_{3,1}+15\,u_{3,2}^2.
\]
If $u_{3,2}\neq0$ then $u_{2,2}=0$ simplifies the equation as
\[
12\,u_{3,2}^2=15\,u_{3,2}^2,
\]
so $u_{3,2}=0$.

By checking $u_{6,3}$ and using $u_{3,2}=0$, we get
\[
u_{6,3}=\frac{945}{4}\,u_{2, 2}^3=225\,u_{2,2}^3,
\]
so $u_{2,2}=0$. We get a form
\[
u=xy+\frac{1}{6}x^3+\frac{I^{\sf rel}_{3,1}}{6}x^3 y+O_x(4)+O_y(3).
\]
The stabilizer, as mentioned in the previous section, is
\[
x'=\mu x, \ \ y'=\mu^2 y, \ \ u'=\mu^3 u,\ \ \mu\neq0.
\]
So when $I^{\sf rel}_{3,1}\neq0$ we can normalize it to 1 and get N7 in \cite{Eastwood-Ezhov-1999}
\[
u=xy+\frac{1}{6}x^3+\frac{1}{6}x^3 y+O_x(4)+O_y(3).
\]
When $I^{\sf rel}_{3,1}\equiv0$ we get N8 in \cite{Eastwood-Ezhov-1999}
\[
u=xy+\frac{1}{6}x^3+O_x(4)+O_y(3),
\]
which is claimed to be the Cayley surface $u=xy+x^3$.

\section{Vanishing $I^{\sf rel}_{0,3}\equiv0\equiv I^{\sf rel}_{3,0}$, $B_3$}\label{sect-v03-v30}
In this branch $I^{\sf rel}_{0,3}\equiv0\equiv I^{\sf rel}_{3,0}$. There are only finitely many independent jet coordinates $u_{\leqslant 2,\leqslant2}$. Among them we have normalized
\[
\aligned
u_{0,2}&=0,  & u_{1,2}&=0,  & u_{2,2}& ,\\
u_{0,1}&=0,  & u_{1,1}&=1,  & u_{2,1}&=0,\\
u_{0,0}&=0,  & u_{1,1}&=0,  & u_{2,0}&=0,
\endaligned
\]
except $u_{2,2}$. The form
\[
u=xy+\frac{F_{2,2}}{4}x^2 y^2+O_x(3)+O_y(3),
\]
has stabilizer
\[
x'=\mu x, \ \ y'=\lambda y, \ \ u'=\mu\lambda u, \ \ \mu,\lambda\in\bc^*.
\]
Thus $F_{2,2}$ is a relative invariant. We denote $F_{2,2}$ by $I^{\sf rel}_{2,2}$. When $I^{\sf rel}_{2,2}\neq0$ we can normalize it to 1 and get
\[
u=xy+\frac{1}{4}x^2 y^2+O_x(3)+O_y(3).
\]
Analysing $I^{\sf rel}_{0,3}\equiv I^{\sf rel}_{3,0}$ and the recurrence relations, for homogeneous models, we get
\[
u=xy+\frac{1}{4}x^2 y^2+O(5)=2-2\sqrt{1-xy},
\]
which is N9 in \cite{Eastwood-Ezhov-1999}.

When $I^{\sf rel}_{2,2}\equiv 0$ we may verify that
\[
u=xy+O_x(3)+O_y(3)=xy+O(5)=xy
\]
which is N10 in \cite{Eastwood-Ezhov-1999}.

\section{Conclusion}
We discover all models in \cite{Eastwood-Ezhov-1999}.
\[
\xymatrix{
&&
{\tiny G_{3,0}\neq0\neq G_{0,3}~\boxed{N1,N2,N3,N4}}
&&
{\tiny G_{4,0}\neq0~\boxed{N5,N6}}
&&
{\tiny G_{3,1}\neq0~\boxed{N7}}
&&
\\
&&
&&
&&
&&
\\
\ar[uurr]^{B_1}
\ar[rr]^{B_2\,\,\,\,\,\,\,\,\,\,\,\,}
{\tiny \boxed{\sf root}}
\ar[ddrr]^{B_3}
&&
\ar[uurr]_{B_{2\cdot 1}}
\ar[rr]_{B_{2\cdot 2}}
G_{3,0}\neq0\equiv G_{0,3}
&&
\ar[uurr]_{B_{2\cdot2\cdot 1}}
\ar[rr]_{B_{2\cdot2\cdot2}}
G_{4,0}\equiv 0
&&
{\tiny G_{3,1}\,\equiv\,0~\boxed{N8}}
&&
\\
&&
&&
&&
&&
\\
&&
G_{3,0}\equiv 0\equiv G_{0,3}
\ar[rr]^{B_{3\cdot 1}}
\ar[rrdd]^{B_{3\cdot 2}}
&&
{\tiny G_{2,2}\neq 0~\boxed{N9}}
&&
&&
\\
&&
&&
\\
&&
&&
{\tiny G_{2,2}\equiv 0~\boxed{N10}}
}
\]

\end{document}